\documentclass[a4paper, 11pt]{article}
\usepackage[left=2cm, right=2cm, top=3cm, bottom=3cm]{geometry}
\usepackage{amssymb, amsmath, mathbbol, amsfonts}
\usepackage{mathtools}
\usepackage{latexsym}
\usepackage{mathrsfs}
\usepackage{bbm}
\usepackage{amsthm}
\usepackage{titleref}
\usepackage{prettyref}
\usepackage{graphicx}
\usepackage{enumerate}
\usepackage{lmodern}
\usepackage{trfsigns}
\usepackage{tikz}
\usetikzlibrary{decorations,arrows}
\usetikzlibrary{decorations.pathmorphing}
\usepgflibrary{decorations.pathreplacing}

\newcommand{\ii}{\mathrm{i}}

\newcommand{\rnzu}{R_{N}(\mathbf{z},\mathfrak{u})}
\newcommand{\sn}{\sum_{\substack{-\frac{N}{2} < h \leq \frac{N}{2} \\ h \neq 0}}{\frac{1}{|h|}}}
\newcommand{\snzo}[1]{\sum\limits_{\substack{-\frac{N}{2} < h \leq \frac{N}{2} \\ h \neq 0}}{\frac{\mathrm{e}^{2 \pi \ii h k b^{w_{#1{}}}z/N}}{|h|}}}
\newcommand{\snz}[1]{\sum\limits_{\substack{-\frac{N}{2} < h \leq \frac{N}{2} \\ h \neq 0}}{\frac{\mathrm{e}^{2 \pi \ii h k b^{w_{#1{}}}z_{#1{}}/N}}{|h|}}}
\newcommand{\snzabs}[1]{\sum\limits_{\substack{-\frac{N}{2} < h \leq \frac{N}{2} \\ h \neq 0}}{\frac{\left|\mathrm{e}^{2 \pi \ii h k b^{w_{#1{}}}z_{#1{}}/N}\right|}{|h|}}}
\newcommand{\snb}[1]{\sum_{\substack{-\frac{N}{2} < h \leq \frac{N}{2} \\ h \neq 0}}{\frac{\mathrm{e}^{2 \pi \ii h k b^{#1{}}/N}}{|h|}}}
\newcommand{\snbn}[1]{\sum_{\substack{-\frac{N}{2} < h \leq \frac{N}{2} \\ h \neq 0}}{\frac{\mathrm{e}^{2 \pi \ii h k b^{#1{}}}}{|h|}}}
\newcommand{\dis}[1]{\mathrm{discr}(#1{},P_N(\mathbf{z}))}
\newcommand{\bsgamma}{\boldsymbol{\gamma}}
\newcommand{\bseta}{\boldsymbol{\eta}}

\newtheorem{thm}{Theorem}[section]
\newtheorem{cor}[thm]{Corollary}

\newtheorem{lem}[thm]{Lemma}
\newtheorem{defi}{Definition}
\newtheorem{alg}{Algorithm}
\newtheorem{rem}{Remark}

\newenvironment{proof1}{\begin{trivlist}\item[\hskip\labelsep{\it Proof.}]}{$\hfill\Box$\end{trivlist}}

\numberwithin{equation}{section}
\allowdisplaybreaks

\title{A reduced fast component-by-component construction of lattice point sets with small weighted star discrepancy}
\author{Ralph Kritzinger\footnote{R. Kritzinger is supported by the Austrian Science Fund (FWF): Project F5509-N26, which is a part of the Special Research Program ``Quasi-Monte Carlo Methods: Theory and Applications''.} \ and Helene Laimer\footnote{H. Laimer is supported by the Austrian Science Fund (FWF): Project F5506-N26, which is a part of the Special Research Program ``Quasi-Monte Carlo Methods: Theory and Applications''.}}
\date{}

\begin{document}
\allowdisplaybreaks
\setcounter{page}{1}
\newpage
\pagestyle{plain}

\maketitle

\begin{abstract}
\noindent
The weighted star discrepancy of point sets appears in the weighted Koksma-Hlawka inequality and thus is a measure for the quality of point sets with respect to their performance in quasi-Monte Carlo algorithms. A special choice of point sets are lattice point sets whose generating vector can be obtained one component at a time such that the resulting lattice point set has a small weighted star discrepancy.

In this paper we consider a reduced fast component-by-component algorithm which significantly reduces the construction cost for such generating vectors provided that the weights decrease fast enough.
\end{abstract}


\section{Introduction}
Given an $N$-element multiset of points $\{ \mathbf{x}_0, \dotsc , \mathbf{x}_{N-1} \} \in [0,1)$, we may approximate integrals over the $s$-dimensional unit cube by a quasi-Monte Carlo (QMC) rule, i.e.,
$$\int_{[0,1]^s}{f(\mathbf{x})} \mathrm{d} \, \mathbf{x} \approx \frac{1}{N} \sum_{n=0}^{N-1}{f(\mathbf{x}_n)}.$$
For detailed information on QMC-integration see \cite{DP, L, LP, Nied2}.

In 1998 Sloan and Wo\'{z}niakowski \cite{SW} introduced the concept of weighted function spaces where each group of coordinates is equipped with some weight according to its importance. Denote the set $\{ 1, \dotsc , s \}$ by $[s]$ and let $\bsgamma = (\bsgamma_{\mathfrak{u}})_{\mathfrak{u} \subseteq [s]}$ be a weight sequence of non-negative real numbers, which model the importance of the projection of the integrands $f$ in the weighted function space onto the variables $x_j$ for $j \in \mathfrak{u}$. A small weight $\bsgamma_{\mathfrak{u}}$ means that the projection onto the variables in $\mathfrak{u}$ contributes little to the integration problem. In the present work we consider a special choice of weights, so-called product weights $(\gamma_j)_{j \geq 1}$, where $\bsgamma_{\mathfrak{u}} = \prod_{j \in \mathfrak{u}}{\gamma_j}$ and $\bsgamma_{\emptyset} := 1$, and in particular, the weight $\gamma_j$ is associated with the variable $x_j$.

In this paper we assume that $\bsgamma = (\gamma_j)_{j \geq 1}$ is a non-increasing sequence of positive weights with $\gamma_j \leq 1$ and $(\bsgamma_{\mathfrak{u}})_{\mathfrak{u} \subseteq [s]}$ are the corresponding product weights. Such weights are useful when considering functions whose dependence on successive variables is decreasing.\\[.5cm]
A particularly important kind of point sets for QMC-integration are so-called lattice point sets. They originated independently from Hlawka \cite{Hlawka2} and Korobov \cite{K}. A lattice point set $P_N(\mathbf{z}) = \{ \mathbf{x}_0, \dotsc , \mathbf{x}_{N-1} \}$ can be constructed with the aid of a generating vector $\mathbf{z}$. For a positive integer $N \geq 2$ and a vector $\mathbf{z} \in \{ 1, \dotsc , N-1\}^s$ the corresponding lattice point set is of the form
$$P_N(\mathbf{z}) = \left\{ \left\{ \frac{k}{N} \mathbf{z} \right\} : k = 1 , \dotsc , N-1\right\}.$$
The brackets $\{ . \}$ around $\frac{k}{N} \mathbf{z}$ indicate that we take the fractional part of each point. For vectors, $\{ . \}$ is applied component-wise. See \cite{LP, Nied2, SloanJ}.

We want to measure the quality of lattice point sets $P_N(\mathbf{z})$ with respect to their performance in a QMC rule. Therefore we define the weighted star discrepancy.
\begin{defi}\label{StarDisc}
Let $\bsgamma = (\bsgamma_{\mathfrak{u}})_{\mathfrak{u} \subseteq [s]}$ be a weight sequence and let $P_N = \left\{ \mathbf{x}_0 , \dotsc , \mathbf{x}_{N-1} \right\} \subseteq [0,1]^s$ be an $N$-element point set. The local discrepancy of the point set $P_N$ at $\mathbf{x} = (x_1, \dotsc, x_s) \in [0,1]^s$ is defined as
$$\mathrm{discr}(\mathbf{x},P_N(\mathbf{z})) := \frac{1}{N}\sum\limits_{\mathbf{p} \in P_N}{\chi_{[\mathbf{0},\mathbf{x})}(\mathbf{p})} - \prod_{j=1}^s{x_j},$$
where $\chi_{[\mathbf{0},\mathbf{x})}$ denotes the characteristic function of $[\mathbf{0},\mathbf{x})$. The weighted star discrepancy of $P_N$ is then defined as
\begin{equation*}
D_{N, \bsgamma}^*(P_N) := \sup_{\mathbf{x} \in (0,1]^s}{ \max_{\emptyset \neq \mathfrak{u} \subseteq [s]}{\bsgamma_{\mathfrak{u}}|\mathrm{discr}((\mathbf{x}_{\mathfrak{u}},\mathbf{1}),P_N)|}}.
\end{equation*}
\end{defi}
We denote the weighted star discrepancy of a lattice point set corresponding to some generating vector $\mathbf{z}$ by $D_{N, \bsgamma}^*(\mathbf{z})$, as this $P_N(\mathbf{z})$ is completely determined by $\mathbf{z}$. To see why the weighted star discrepancy is a measure for the quality of our point sets we study the following identity of Hlawka \cite{Hlawka1} and Zaremba \cite{Z} (see also \cite{DP, LP}), given by 
\begin{align*}
Q_{N,s}{(f)} - I_s(f) = \sum_{\emptyset \neq \mathfrak{u} \subseteq [s]}{(-1)^{|\mathfrak{u}|} \bsgamma_{\mathfrak{u}}\int_{[0,1]^{|\mathfrak{u}|}}{\dis{(\mathbf{x}_{\mathfrak{u}},\mathbf{1})}\bsgamma_{\mathfrak{u}}\frac{\partial^{|\mathfrak{u}|}}{\partial\mathbf{x}_{\mathfrak{u}}}f\left(\mathbf{x}_{\mathfrak{u}},\mathbf{1}\right) \, \mathrm{d}\,\mathbf{x}_{\mathfrak{u}}}},
\end{align*}
where $Q_{N,s}(f) = \frac{1}{N}\sum_{j=1}^s{f(\mathbf{x}_j)}$ denotes the QMC-rule, $I_s = \int_{[0,1]^s}{f(\mathbf{x})} \, \mathrm{d}\mathbf{x}$ the integral operator and $(\mathbf{x}_{\mathfrak{u}},\mathbf{1})$ the vector $(\tilde{x}_1, \dotsc , \tilde{x}_s)$ with $\tilde{x}_j = x_j$ if $j \in \mathfrak{u}$ and $\tilde{x}_j = 1$ if $j \notin \mathfrak{u}$.

Applying Hölder's inequality as in \cite{DP, SW} for integrals and sums we obtain
\begin{align}\label{eq:weightedKoksmaHlawka}
|Q_{N,s}{(f)} - I_s(f)| \leq D_{N, \bsgamma}^*(\mathbf{z}) \|f\|_{\bsgamma},
\end{align}
where $\|.\|_{\bsgamma}$ is some norm dependent on $\bsgamma$ but independent of the point set $P_N{(\mathbf{z})}$. If $f$ is sufficiently smooth $\|f\|_{\bsgamma}$ coincides with the weighted variation of $f$ in the sense of Hardy and Krause. The first factor in \eqref{eq:weightedKoksmaHlawka} is the weighted star discrepancy of the point set $P_N(\mathbf{z})$ and depends only on $P_N(\mathbf{z})$ and the weights. Thus we see that the quality of a lattice point set $P_N(\mathbf{z})$ is the better the smaller its weighted star discrepancy $D_{N, \bsgamma}^*(\mathbf{z})$. We want to find lattice point sets $P_N(\mathbf{z})$ with small weighted star discrepancy.

As no explicit constructions for good lattice point sets are known for dimensions $s > 2$, one usually employs computer search algorithms to find good generating vectors. There exist many papers on the construction of generating vectors for lattice point sets with a small weighted star discrepancy: Joe \cite{J} has given a component-by-component construction for generating vectors of lattice point sets with a prime number $N$ of points, which have a weighted star discrepancy of order $N^{-1+\delta}$ for any $\delta > 0$. Their generating vector has a construction cost of order $s N \log N$, where an approach of Nuyens and Cools \cite{NC2} can be used to reduce the construction cost.

In \cite{SJ1} Joe and Sinescu have achieved the same results for a composite number of lattice points and product weights. Finally in \cite{SJ2} they considered general weights and a prime number of points.

Dick et al. \cite{DKLP} have given a reduced fast algorithm for the construction of generating vectors of lattice point sets with $N$ a prime power. They varied the size of the search space for each coordinate according to its importance and considered the worst-case error of integration in a Korobov space to measure the quality of their lattice point sets.

Let $b$ be an arbitrary prime number and $m$ a positive integer. In the present work we consider lattice point sets with $N = b^m$ elements and study their weighted star discrepancy. As mentioned before, the generating vector $\mathbf{z} = (z_1, \dotsc , z_s)$ of such lattice point sets can be obtained one component at a time. When using the standard component-by-component construction, in the following frequently abbreviated by CBC construction, each component is chosen from $\left\{ z \in \left\{ 1, 2, \dotsc , b^m - 1 \right\} : \gcd{(z,b^m)} = 1 \right\}$. As done in \cite{DKLP} for the worst-case error, we speed up the construction of such generating vectors by reducing the search space for each component, while still achieving a small weighted star discrepancy of the corresponding lattice rule. To this end we define non-decreasing $0 \leq w_1 \leq w_2 \leq \dotsc \in \mathbbm{N}$ and set
\begin{equation*}
\mathcal{Z}_{N,w_j} :=	\begin{cases}
									\left\{ z \in \left\{ 1, 2, \dotsc , b^{m-w_j}-1 \right\} : \gcd{(z,b^m)} = 1 \right\}  	& \text{ if } w_j < m, \\
									\{ 1\} 																																									& \text{ if } w_j \geq m.
												\end{cases}
\end{equation*}
Note that these sets have cardinality $b^{m-w_j-1}(b-1)$, for $w_j < m$. In what follows we denote by $\mathcal{Z}_{N,\mathbf{w}}^s$ the cartesian product $b^{w_1} \mathcal{Z}_{N,w_1} \times \dotsc \times b^{w_s} \mathcal{Z}_{N,w_s}$, where $b^{w_j} \mathcal{Z}_{N,w_j}$ means that every element of $\mathcal{Z}_{N,w_j}$ is multiplied by $b^{w_j}$. We denote by $\mathbf{z} \in \mathcal{Z}_{N,\mathbf{w}}^s$ a vector $\mathbf{z} = (b^{w_1} z_1, \dotsc , b^{w_s} z_s)$, with $z_j \in \mathcal{Z}_{N,w_j}$ for $j \in [s]$. We study the weighted star discrepancy of lattice point sets $P_N(\mathbf{z})$ with generating vectors $\mathbf{z} \in \mathcal{Z}_{N,\mathbf{w}}^s$. Dick et al. \cite{DKLP} have considered the worst-case error for approximating the integral of functions in suitable spaces by a QMC rule based on lattice point sets. Here, in contrast, we study the weighted star discrepancy of these lattice point sets which is another important quality measure. We will see that for sufficiently fast decreasing weights we can construct lattice point sets with small weighted star discrepancy, while significantly reducing the construction cost in comparison to the standard CBC construction.\\[.5cm]
It follows from \cite[Theorem~3.10 and Theorem~5.6]{Nied2} that 
\begin{equation}\label{eq:estStar1}
D_{N, \bsgamma}^*(\mathbf{z}) \leq \sum_{\mathfrak{u} \subseteq [s]}{\gamma_{\mathfrak{u}}\left( 1- \left( 1 - \frac{1}{N} \right) \right)^{|\mathfrak{u}|}} + \frac{1}{2} R_{N,\bsgamma}^s(\mathbf{z}),
\end{equation}
where
\begin{align}\label{eq:endSquared}
R_{N,\bsgamma}^s(\mathbf{z}) = \sum_{\mathfrak{u} \subseteq [s]}{\gamma_{\mathfrak{u}} \rnzu}
\end{align}
and
\begin{equation}\label{rnzu}
\rnzu = \frac{1}{N}\sum_{k=0}^{N-1}{\prod_{j \in \mathfrak{u}}{\left( 1 + \snz{j} \right)}} - 1.
\end{equation}
Using this estimate for the weighted star discrepancy we derive the results in Sections \ref{sectionArihmeticMean}, \ref{sectionReducedCBC} and \ref{sectionReducedFastCBC}.

Finally, we introduce the concept of tractability \cite{NW1, NW2, NW3}. To this end we define the information complexity (often refered to as inverse of the weighted star discrepancy) as
$$N^*(\varepsilon,s) = \min\{N \in \mathbbm{N}_0 : D_{N, \bsgamma}^*(\mathbf{z}) \leq \varepsilon \},$$
which means that $N^*(\varepsilon,s)$ is the minimal number of points required to achieve a weighted star discrepancy of at most $\varepsilon$. Of course we want the information complexitiy to be as small as possible. Therefore we are interested in how fast it increases when $\varepsilon^{-1}$ and $s$ grow. We define the following notions of tractability. We speak of
\begin{itemize}
\item polynomial tractability, if there exist constants $C, \tau_1 > 0$ and $\tau_2 \geq 0$ such that
$$N^*(\varepsilon, s) \leq C \varepsilon^{-\tau_1}s^{\tau_2} \text{ for all } \varepsilon \in (0,1) \text{ and all } s \in \mathbbm{N} \text{ and of }$$
\item strong polynomial tractability, if there exist positive constants $C, \tau$ such that
$$n(\varepsilon, s) \leq C \varepsilon^{-\tau} \text{ for all } \varepsilon \in (0,1) \text{ and all } s \in \mathbbm{N}.$$
\end{itemize}
Roughly speaking, a problem is considered tractable if its information complexity's dependence on $\varepsilon^{-1}$ and $s$ is not exponential.
We will show that the above mentioned reduced fast component-by-component construction finds a generating vector $\mathbf{z}$ of a lattice point set that achieves strong polynomial tractability if
$$\sum_{j=1}^{\infty}{\gamma_j b^{w_j}} < \infty$$
with a construction cost of
$$O\left( N \log N + \min\{s,t\} N + N \sum_{d=1}^{\min\{s,t\}}{(m-w_d)b^{-w_d}} \right)$$
operations, where $t = \max\{j \in \mathbbm{N} : w_{j} < m\}$.

The structure of this paper is as follows. In the next section we derive an upper bound for the arithmetic mean of the weighted star discrepancy over all possible lattice point sets constructed by a generating vector $\mathbf{z} \in \mathcal{Z}_{N,\mathbf{w}}^s$. In Sections \ref{sectionReducedCBC} and \ref{sectionReducedFastCBC} we present a reduced fast CBC construction for generating vectors of lattice point sets with small weighted star discrepancy. Finally, in Section \ref{sectionTractability} we study conditions on the weights $\gamma_j$ and $w_j$ for achieving strong polynomial tractability.

\section{The arithmetic mean over all $\mathbf{z} \in \mathcal{Z}_{N,\mathbf{w}}^s$}\label{sectionArihmeticMean}
First of all we estimate the arithmetic mean of the weighted star discrepancy over all possible generating vectors $\mathbf{z} = (b^{w_1}z_1, \dotsc , b^{w_s}z_s) \in \mathcal{Z}_{N,\mathbf{w}}^s$, proceeding similarly to \cite{Nied2} and \cite{SJ1}. This yields the existence of a lattice point set with small weighted star discrepancy. The upper bound which we obtain for the arithmetic mean is not the same as for the reduced CBC construction in the next section. Nonetheless, we need large parts of the calculation of the present section to obtain the estimate in Section \ref{sectionReducedCBC}.
\begin{thm}\label{thm:boundMean}
Let $N = b^m$, $(w_j)_{j \geq 1}$ and $\in \mathcal{Z}_{N,\mathbf{w}}^s$ be as above and let $m \geq 5$. Then there exists a generating vector $\mathbf{z} = \left( b^{w_1}z_1 , \dotsc , b^{w_s} z_s \right) \in \mathcal{Z}_{N,\mathbf{w}}^s$ whose corresponding lattice rule has weighted star discrepancy
\begin{align*}
D_{N, \bsgamma}^*(\mathbf{z}) &\leq \sum_{\mathfrak{u} \subseteq [s]}{\gamma_{\mathfrak{u}}\left( 1- \left( 1 - \frac{1}{N} \right) \right)^{|\mathfrak{u}|}} + \frac{1}{2} \left(\frac{1}{N} \prod_{j=1}^s{\left( \beta_j + \gamma_j S_N \right)}\right.\\
									& \quad + \frac{1}{N}\sum_{p=0}^{m-1}{b^{m-p-1}(b-1) \prod_{\substack{j=1 \\ w_j \geq m-p}}^{s}{\left( \beta_j + \gamma_jS_N \right)}\prod_{\substack{j=1 \\ w_j < m-p}}^{s}{\beta_j}} - \left.\prod_{j=1}^s{\beta_j} \right),
\end{align*}
with $\beta_j = 1 + \gamma_j$ for all $j \in \mathbbm{N}$ and
\begin{equation}\label{eq:S_N}
S_N = \sn.
\end{equation}. 
\end{thm}
\noindent
\begin{rem} 
Provided that the $\gamma_j$'s are summable the bound in \prettyref{thm:boundMean} is of order $N^{\delta} \log N$ for arbitrary $\delta \in (0,1)$ with an implied constant independent of $N$ and $s$. Furthermore note that if all weights $w_j = 0$ then we obtain the result in \cite[Theorem~1 and Corollary~1]{SJ1}.
\end{rem}
\begin{proof1}
As the first sum in \eqref{eq:estStar1} is independent of $\mathbf{z}$, it is obviously enough to consider the mean
\begin{equation}\label{eq:mean}
M_{N,s,\bsgamma} := \frac{1}{|\mathcal{Z}_{N,\mathbf{w}}^s|}\sum_{\mathbf{z} \in \mathcal{Z}_{N,\mathbf{w}}^s}{R_{N,\bsgamma}^s(\mathbf{z})}
\end{equation}
of the second sum.

We have from \cite[p.~186, eq.~9]{J}
\begin{equation}\label{eq:endSquared1}
\begin{split}
R_{N,\bsgamma}^s(\mathbf{z}) &= \frac{1}{N}\sum_{k=0}^{N-1}{\prod_{j=1}^s{\left( \beta_j + \gamma_j \snz{j} \right)}} - \prod_{j=1}^s{\beta_j}\\
											&= \frac{1}{N}\prod_{j=1}^s{\left( \beta_j + \gamma_j S_N \right)} + \frac{1}{N}\sum_{k=1}^{N-1}{\prod_{j=1}^s{\left( \beta_j + \gamma_j \snz{j} \right)}} - \prod_{j=1}^s{\beta_j}.
\end{split}
\end{equation}
Thus
\begin{equation*}
\begin{split}
M_{N,s,\bsgamma} &= \frac{1}{N}\prod_{j=1}^s{\left( \beta_j + \gamma_j S_N \right)} + \frac{1}{N}\sum_{k=1}^{N-1}{\prod_{j=1}^s{\left(\frac{1}{|\mathcal{Z}_{N, w_j}|}\sum_{z_j \in \mathcal{Z}_{N, w_j}}{\left( \beta_j + \gamma_j \snz{j} \right)}\right)}} - \prod_{j=1}^s{\beta_j} \\
												&= \frac{1}{N}\prod_{j=1}^s{\left( \beta_j + \gamma_j S_N \right)} + \frac{1}{N}\sum_{k=1}^{N-1}{\prod_{j=1}^s{\left( \beta_j + \frac{\gamma_j}{|\mathcal{Z}_{N, w_j}|}\sum_{z_j \in \mathcal{Z}_{N, w_j}}{ \snz{j} }\right)}} - \prod_{j=1}^s{\beta_j}.
\end{split}
\end{equation*}
To avoid lengthy formulas we use the following abbreviations:
\begin{equation}\label{eq:tnwjk}
T_{N,w_j}(k) := \sum_{z_j \in \mathcal{Z}_{N, w_j}}{ \snz{j} }
\end{equation}
and
\begin{equation}\label{eq:lnsgamma}
L_{N,s,\bsgamma} := \frac{1}{N}\sum_{k=1}^{N-1}{\prod_{j=1}^s{\left( \beta_j + \frac{\gamma_j}{|\mathcal{Z}_{N, w_j}|}T_{N,w_j}(k) \right)}}.
\end{equation}
Then we have
\begin{equation}\label{eq:meanKurz}
M_{N,s,\bsgamma} = \frac{1}{N}\prod_{j=1}^s{\left( \beta_j + \gamma_j S_N \right)} + L_{N,s,\bsgamma} - \prod_{j=1}^s{\beta_j}.
\end{equation}
We study $T_{N,w_j}(k)$ distinguishing the two cases $w_j \geq m$ and $w_j < m$.\\[.5cm]
\noindent
{\bfseries Case 1:} $w_j \geq m$. This yields $\mathcal{Z}_{N, w_j} = \{1\}$ and thus
\begin{equation}\label{eq:geqm}
T_{N,w_j}(k) = \snb{w_j} = \snbn{w_j-m} = \sn = S_N.
\end{equation}
{\bfseries Case 2:} $w_j < m$. Then $\mathcal{Z}_{N, w_j} = \left\{ z \in \left\{ 1, 2, \dotsc , b^{m-w_j}-1 \right\} : \gcd{(z,N)} = 1 \right\}$.
According to \eqref{eq:lnsgamma} we have to calculate $T_{N,w_j}(k)$ only for $k \in \{ 1, \dotsc , b^m - 1 \}$. We display these $k$ as $k = q b^{m-w_j} + r$ with $q \in \{ 0, \dotsc , b^{w_j}-1 \}$, $r \in \{0, \dotsc, b^{m-w_j}-1\}$ and $(q,r) \neq (0,0)$. Then
\begin{equation}\label{eq:<m}
\begin{split}
T_{N,w_j}(k)	&= \sn  \sum_{z_j \in \mathcal{Z}_{N, w_j}}{\mathrm{e}^{2 \pi \ii h (q b^{m-w_j} + r) b^{w_j} z_j/N}} = \sn  \sum_{z_j \in \mathcal{Z}_{N, w_j}}{\mathrm{e}^{2 \pi \ii h q z_j} \mathrm{e}^{2 \pi \ii h r z_j/b^{m-w_j}}} \\
							&= \sn  \sum_{z_j \in \mathcal{Z}_{N, w_j}}{ \mathrm{e}^{2 \pi \ii h r z_j/b^{m-w_j}}}.
\end{split}
\end{equation}

If $r=0$, i.e. $k$ a multiple of $b^{m-w_j}$, this yields
\begin{equation}\label{eq:r=0}
T_{N,w_j}(k) = \sn\sum_{z_j \in \mathcal{Z}_{N, w_j}}{1} = |\mathcal{Z}_{N, w_j}| S_N.
\end{equation}

Next we investigate $r \in \{ 1, \dotsc , b^{m-w_j}-1 \}$. For any $z_j \in \{0, \dotsc , b^{m-w_j}-1\}$ we find $\gcd{(z_j,N)} = \gcd{(z_j,b^{m-w_j})} \in \left\{ 1, b, b^2, \dotsc , b^{m-w_j-1} \right\}$ and hence
\begin{equation*}
\sum_{d | \gcd{(z_j,N)}}{\mu(d)} = \sum_{d | \gcd{(z_j,b^{m-w_j})}}{\mu(d)} = \begin{cases}
																																								1 & \text{iff } \gcd{(z_j,N)} = \gcd{(z_j,b^{m-w_j})} = 1,\\
																																								0 & \text{otherwise,}
																																							\end{cases}
\end{equation*}
where $\mu$ denotes the Möbius function.

For any $z_j \in \{ 1, \dotsc , b^{m-w_j} - 1 \}$ this implies $z_j \in \mathcal{Z}_{N, w_j}$ if and only if $\sum\limits_{d | \gcd{(z_j,b^{m-w_j})}}{\mu(d)} = 1$. Inserting this fact into \eqref{eq:<m} we have
\begin{equation}\label{eq:<m2}
T_{N,w_j}(k) = \sn  \sum_{z_j = 1}^{b^{m-w_j}-1}{\mathrm{e}^{2 \pi \ii h r z_j/b^{m-w_j}}\sum_{d | \gcd{(z_j,b^{m-w_j})}}{\mu(d)}}.
\end{equation}
Studying the two inner sums we find
\begin{equation}\label{eq:r>0}
\begin{split}
\sum_{z_j = 1}^{b^{m-w_j}-1}{\mathrm{e}^{2 \pi \ii h r z_j/b^{m-w_j}}\sum_{d | \gcd{(z_j,b^{m-w_j})}}{\mu(d)}}	&= \sum_{d | b^{m-w_j}}{\mu(d)\sum_{\substack{z_j = 1 \\ d | z_j}}^{b^{m-w_j}-1}{\mathrm{e}^{2 \pi \ii h r z_j/b^{m-w_j}}}} \\
																																																												&= \sum_{d | b^{m-w_j}}{\mu(d)\sum_{a = 1}^{\frac{b^{m-w_j}}{d}}{\mathrm{e}^{2 \pi \ii h r a d/b^{m-w_j}}}},
\end{split}
\end{equation}
where the latter equality holds since $a \in \left\{ 1, \dotsc , \frac{b^{m-w_j}}{d} \right\}$ yields $a d \in \left\{ d, 2 d, \dotsc , b^{m-w_j} \right\}$\newline
$= \left\{ 1 \leq z_j \leq b^{m-w_j} - 1 : d | z_j \right\} \cup \left\{ b^{m-w_j} \right\}$ and 
\begin{align*}
\sum_{d | b^{m-w_j}}{\mu(d)} = 0,
\end{align*}
since $w_j < m.$

Changing the order of summation we obtain with \eqref{eq:r>0}
\begin{equation*}
\sum_{z_j = 1}^{b^{m-w_j}-1}{\mathrm{e}^{2 \pi \ii h r z_j/b^{m-w_j}}\sum_{d | \gcd{(z_j,b^{m-w_j})}}{\mu(d)}} = \sum_{d | b^{m-w_j}}{\mu\! \left(\frac{b^{m-w_j}}{d}\right)\sum_{a = 1}^{d}{\mathrm{e}^{2 \pi \ii h r a/d}}} = \sum_{\substack{d | b^{m-w_j} \\ d | h r}}{d \, \mu\! \left(\frac{b^{m-w_j}}{d}\right)}.
\end{equation*}
With \eqref{eq:<m2} this leads to
\begin{equation*}
T_{N,w_j}(k) = \sn  \sum_{\substack{d | b^{m-w_j} \\ d | h r}}{d \, \mu\! \left(\frac{b^{m-w_j}}{d}\right)} = \sum_{d | b^{m-w_j}}{d \, \mu\! \left(\frac{b^{m-w_j}}{d}\right)} \sum_{\substack{-\frac{N}{2} < h \leq \frac{N}{2} \\ h \neq 0 \\ d | h r}}{\frac{1}{|h|}}.
\end{equation*}
Using that $d | h r$ is equivalent to $\frac{d}{\gcd{(d,r)}} | h$ we display $T_{N,w_j}(k)$ as
\begin{equation}\label{eq:<m3}
T_{N,w_j}(k) = \sum_{d | b^{m-w_j}}{d \, \mu\! \left(\frac{b^{m-w_j}}{d}\right)} \sum_{\substack{-\frac{N}{2} < h \leq \frac{N}{2} \\ h \neq 0 \\ \frac{d}{\gcd{(d,r)}} | h}}{\frac{1}{|h|}}.
\end{equation}

For further investigation of $T_{N,w_j}(k)$ we first study sums of the same type as the inner sum in \eqref{eq:<m3}. For any positive integer $a$ we have 
\begin{equation}\label{eq:NiedSn}
\sum_{\substack{-\frac{N}{2} < h \leq \frac{N}{2} \\ h \neq 0 \\ a | h}}{\frac{1}{|h|}} = \sum_{\substack{-\frac{N}{2} < a p \leq \frac{N}{2} \\ p \neq 0}}{\frac{1}{a |p|}} = \frac{1}a{} \sum_{\substack{-\frac{N}{2a} < p \leq \frac{N}{2a} \\ p \neq 0}}{\frac{1}{|p|}} = \frac{1}{a} S_{\frac{N}{a}}, 
\end{equation}
where $S_{\frac{N}{a}}$ is defined analogously to \eqref{eq:S_N}. Combining \eqref{eq:NiedSn} with \eqref{eq:<m3} we obtain
\begin{equation}\label{eq:<m4}
\begin{split}
T_{N,w_j}(k)	&= \sum_{d | b^{m-w_j}}{d \, \mu\! \left(\frac{b^{m-w_j}}{d}\right) \frac{\gcd{(d,r)}}{d}S_{\frac{N}{d}\gcd{(d,r)}}} = \sum_{d | b^{m-w_j}}{ \mu\! \left(\frac{b^{m-w_j}}{d}\right) \gcd{(d,r)}S_{\frac{N}{d}\gcd{(d,r)}}}\\
							&= \gcd{(b^{m-w_j},r)} S_{b^{w_j}\gcd{(b^{m-w_j},r)}} - \gcd{(b^{m-w_j-1},r)} S_{b^{w_j+1}\gcd{(b^{m-w_j-1},r)}} \\
							&= b^{\nu}(S_{b^{w_j + \nu}} - S_{b^{w_j + \nu + 1}}),
\end{split}
\end{equation}
with $\nu \in \{ 0, \dotsc , m-w_j-1\}$.

Summarizing, we have for $k \in \left\{ 1, \dotsc , b^{m}-1 \right\}$
\begin{equation}\label{eq:tnwjk1}
T_{N,w_j}(k) =	\begin{cases}
									S_N																																																		& \text{ if } \quad w_j \geq m,\\
									|\mathcal{Z}_{N,w_j}|S_N 																																							& \text{ if } \quad w_j < m \text{ and } k \equiv 0 \ ( \mathrm{mod} \, b^{m-w_j}), \\
									b^{\nu}(S_{b^{w_j + \nu}} - S_{b^{w_j + \nu + 1}}) 																										&\\
									\qquad \text{ with } b^{\nu} = \gcd{(b^{m-w_j},r)}																										& \text{ if } \quad w_j < m \text{ and } k \not\equiv 0 \ ( \mathrm{mod} \, b^{m-w_j}).
								\end{cases}
\end{equation}

Let us choose $t \in \mathbbm{N}_0$ such that $w_j < m$ for all $j \leq t$ and $w_{t+1} \geq m$. (If $t = 0$, then $w_j \geq m$ for all $j \in \mathbbm{N}$. In that case we obtain the generating vector $\mathbf{z} = (b^{w_1}, \dotsc , b^{w_s})$.) With this we are able to write $L_{N,s,\bsgamma}$ from formula \eqref{eq:lnsgamma} as
\begin{equation}\label{eq:lnsgamma1}
\begin{split}
L_{N,s,\bsgamma}	&= \frac{1}{N}\sum_{k=1}^{N-1}{\prod_{j=1}^{\min\{t, s\}}{\left( \beta_j + \frac{\gamma_j}{|\mathcal{Z}_{N, w_j}|}T_{N,w_j}(k) \right)} \prod_{j=t+1}^s{\left( \beta_j + \frac{\gamma_j}{|\mathcal{Z}_{N, w_j}|}T_{N,w_j}(k) \right)}} \\
									&= \frac{1}{N}\prod_{j=t+1}^s{\left( \beta_j + \gamma_j S_N \right)}\sum_{k=1}^{N-1}{\prod_{j=1}^{\min\{t, s\}}{\left( \beta_j + \frac{\gamma_j}{|\mathcal{Z}_{N, w_j}|}T_{N,w_j}(k) \right)} }.
\end{split}
\end{equation}
Next we aim at finding bounds for $\frac{T_{N,w_j}(k)}{|\mathcal{Z}_{N, w_j}|}$ for $w_j < m$.

If $k$ is a multiple of $b^{m-w_j}$ we see immediately from \eqref{eq:tnwjk1} that 
$$\frac{T_{N,w_j}(k)}{|\mathcal{Z}_{N, w_j}|} = \frac{|\mathcal{Z}_{N, w_j}| S_N}{|\mathcal{Z}_{N, w_j}|} = S_N.$$

If $k$ is not a multiple of $b^{m-w_j}$, we use a formula from Niederreiter \cite{Nied1} for $S_n$ with arbitrary $n \in \mathbbm{N}$, given by
\begin{equation}\label{eq:S_m}
S_n = 2 \log{n} + 2 \gamma - \log 4 + \varepsilon{(n)},
\end{equation}
where $\gamma$ denotes the Euler-Mascheroni constant $\gamma = \lim\limits_{l \rightarrow \infty}{\left( \sum\limits_{k=1}^l{\frac{1}{k}} - \log{l} \right)} \approx 0.577216\dotsc$ and 
\begin{equation}\label{eq:epsilon}
\begin{cases}
	-\frac{4}{n^2} < \varepsilon(n) \leq 0,							& \text{ if } n \text{ is even,}\\
	-\frac{3}{n^2} < \varepsilon(n) < \frac{1}{n^2},			& \text{ if } n \text{ is odd.} 
\end{cases}
\end{equation}
From \eqref{eq:tnwjk1} we know
\begin{equation}\label{eq:upperBoundTnwjk}
\begin{split}
T_{N,w_j}(k) = b^{\nu}(S_{b^{w_j + \nu}} - S_{b^{w_j + \nu + 1}}) < 0.
\end{split}
\end{equation}

With $m \geq 5$ we find $-2 < \frac{T_{N,w_j}(k)}{|\mathcal{Z}_{N, w_j}|} < 0$ for $w_j < m$ and $k$ not a multiple of $b^{m-w_j}$ as follows. The upper bound follows immediately from \eqref{eq:upperBoundTnwjk}. It remains to show the lower bound. First we consider $T_{N,w_j}(k)$ using \eqref{eq:S_m}. We have
\begin{equation*}
\begin{split}
T_{N,w_j}(k)	&= b^{\nu}(S_{b^{w_j + \nu}} - S_{b^{w_j + \nu + 1}}) = b^{\nu}\left( -2 \log{b} + \varepsilon{(b^{w_j + \nu})} - \varepsilon{(b^{w_j + \nu + 1})}\right) \\
							&= -2 b^{\nu} \log b + b^{\nu} \left(\varepsilon{(b^{w_j + \nu})} - \varepsilon{(b^{w_j + \nu + 1})}\right).
\end{split}
\end{equation*}
With \eqref{eq:epsilon} we obtain
\begin{equation*}
\begin{split}
\left|b^{\nu} \left(\varepsilon{(b^{w_j + \nu})} - \varepsilon{(b^{w_j + \nu + 1})}\right)\right|	&\leq \left|b^{\nu} \left(\varepsilon{(b^{w_j + \nu})}\right)\right| + \left|b^{\nu} \left(\varepsilon{(b^{w_j + \nu + 1})}\right)\right| \leq 4 b^{-2 w_j - \nu}\left( 1 + \frac{1}{b^2} \right).
\end{split}
\end{equation*}
Thus
\begin{equation*}
\begin{split}
\frac{T_{N,w_j}(k)}{|\mathcal{Z}_{N, w_j}|} &\geq -\frac{b^{w_j-m+1}}{b-1}2 b^{\nu}\log b - \frac{b^{w_j-m+1}}{b-1}4 b^{-2 w_j - \nu}\left( 1 + \frac{1}{b^2} \right).
\end{split}
\end{equation*}
Recall from \eqref{eq:tnwjk1} that $\nu = \log_b{(\gcd{(b^{m-w_j},r)})} \in \left\{0,1, \dotsc , m-w_j-1 \right\}$. Thus
\begin{equation*}
\begin{split}
\frac{T_{N,w_j}(k)}{|\mathcal{Z}_{N, w_j}|} &\geq -2 b^{w_j-m+1+m-w_j-1} \frac{\log b}{b-1} - 4b^{-w_j-m+1-\nu}\frac{1}{b-1}\left( 1 + \frac{1}{b^2} \right) \\
																						&\geq -2 \frac{\log b}{b-1} - 4b^{-m+1}\frac{1}{b-1}\left( 1 + \frac{1}{b^2} \right).
\end{split}
\end{equation*}
Now, with the assumption $m \geq 5$,
\begin{equation*}
\begin{split}
\frac{T_{N,w_j}(k)}{|\mathcal{Z}_{N, w_j}|} &\geq -2 \frac{\log b}{b-1} - 4b^{-5+1}\frac{1}{b-1}\left( 1 + \frac{1}{b^2} \right)\\
																						&\geq -2 \frac{\log 2}{2-1} - 4 \cdot 2^{-5+1}\left( 1 + \frac{1}{2^2} \right) > -2,
\end{split}
\end{equation*}
and hence
$$-2 < \frac{T_{N,w_j}(k)}{|\mathcal{Z}_{N, w_j}|} < 0 \text{ for } w_j < m \text{ and } b^{m-w_j} \nmid k.$$
For any integer $p \in \{ 0, \dotsc , m-1 \}$ with $b^p \mid k$ and $b^{p+1} \nmid k$ the condition $b^{m-w_j} \nmid k$ is equivalent to $m-w_j>p$ or $w_j < m-p$, respectively. Thus we can display \eqref{eq:lnsgamma1} as
\begin{equation*}
\begin{split} 
L_{N,s,\bsgamma}	&= \frac{1}{N}\prod_{j=t+1}^s{\left( \beta_j + \gamma_j S_N \right)}\\
									&\quad \times \sum_{p=0}^{m-1}\sum_{\substack{k=1 \\ b^p \mid k \\ b^{p+1} \nmid \,k }}^{N-1}{ \prod_{\substack{j=1 \\ w_j \geq m-p}}^{\min\{t, s\}}{\left( \beta_j + \frac{\gamma_j}{|\mathcal{Z}_{N, w_j}|}T_{N,w_j}(k) \right)} \prod_{\substack{j=1 \\ w_j < m-p}}^{\min\{t, s\}}{\left( \beta_j + \frac{\gamma_j}{|\mathcal{Z}_{N, w_j}|}T_{N,w_j}(k) \right)}}\\
									&\leq \frac{1}{N}\prod_{j=t+1}^s{\left( \beta_j + \gamma_j S_N \right)}\sum_{p=0}^{m-1}\sum_{\substack{k=1 \\ b^p \mid k \\ b^{p+1} \nmid \,k }}^{N-1}{ \prod_{\substack{j=1 \\ w_j \geq m-p}}^{\min\{t, s\}}{\left( \beta_j + \gamma_jS_N \right)} \prod_{\substack{j=1 \\ w_j < m-p}}^{\min\{t, s\}}{\beta_j}},
\end{split}
\end{equation*}
where the latter estimate holds since $\beta_j > 1$, $-2 < \frac{T_{N,w_j}(k)}{|\mathcal{Z}_{N, w_j}|} < 0$ and $\gamma_j \leq 1$. From 
\begin{equation}\label{eq:cardinality||}
\begin{split}
\left| \left\{ k \in \right. \right. & \left. \left. \left\{ 1, \dotsc , N -1 \right\} : b^p \mid k \text{ and } b^{p+1} \nmid k \right\} \right| \\
																																								& = \left| \left\{ k \in \left\{ 1, \dotsc , b^m-1 \right\} : b^p \mid k \right\} \right| - \left| \left\{ k \in \left\{ 1, \dotsc , b^m-1 \right\} : b^{p+1} \mid k \right\} \right| \\
																																								& = b^{m-p}-1 - \left( b^{m-p-1} - 1 \right) = b^{m-p-1}(b-1)
\end{split}
\end{equation}
we get
\begin{equation*}
L_{N,s,\bsgamma}	\leq \frac{1}{N}\prod_{j=t+1}^s{\left( \beta_j + \gamma_j S_N \right)}\sum_{p=0}^{m-1}{b^{m-p-1}(b-1) \prod_{\substack{j=1 \\ w_j \geq m-p}}^{\min\{t, s\}}{\left( \beta_j + \gamma_jS_N \right)} \prod_{\substack{j=1 \\ w_j < m-p}}^{\min\{t, s\}}{\beta_j}}.
\end{equation*}
Inserting this into \eqref{eq:meanKurz} we obtain for the arithmetic mean
\begin{equation}\label{eq:mean1}
\begin{split}
M_{N,s,\bsgamma}	&= \frac{1}{N}\prod_{j=1}^s{\left( \beta_j + \gamma_j S_N \right)} \\
									&\quad + \frac{1}{N}\prod_{j=t+1}^s{\left( \beta_j + \gamma_j S_N \right)}\sum_{p=0}^{m-1}{b^{m-p-1}(b-1) \prod_{\substack{j=1 \\ w_j \geq m-p}}^{\min\{t, s\}}{\left( \beta_j + \gamma_jS_N \right)} \prod_{\substack{j=1 \\ w_j < m-p}}^{\min\{t, s\}}{\beta_j}} - \prod_{j=1}^s{\beta_j}.
\end{split}
\end{equation}
This proves with \eqref{eq:estStar1} the existence of a vector $\mathbf{z} \in \mathcal{Z}_{N,\mathbf{w}}^s$ such that the weighted star discrepancy $D_{N, \bsgamma}^*(\mathbf{z})$ fulfils
\begin{align}\label{eq:estStar2}
D_{N, \bsgamma}^*(\mathbf{z}) &\leq \sum_{\mathfrak{u} \subseteq [s]}{\gamma_{\mathfrak{u}}\left( 1- \left( 1 - \frac{1}{N} \right) \right)^{|\mathfrak{u}|}} + \frac{1}{2} \left(\frac{1}{N}\prod_{j=1}^s{\left( \beta_j + \gamma_j S_N \right)} \right.\nonumber\\
									& \quad + \frac{1}{N}\prod_{j=t+1}^s{\left( \beta_j + \gamma_j S_N \right)}\sum_{p=0}^{m-1}{b^{m-p-1}(b-1) \prod_{\substack{j=1 \\ w_j \geq m-p}}^{\min\{t, s\}}{\left( \beta_j + \gamma_jS_N \right)} \prod_{\substack{j=1 \\ w_j < m-p}}^{\min\{t, s\}}{\beta_j}} - \left.\prod_{j=1}^s{\beta_j} \right)\nonumber\\
									&\leq \sum_{\mathfrak{u} \subseteq [s]}{\gamma_{\mathfrak{u}}\left( 1- \left( 1 - \frac{1}{N} \right) \right)^{|\mathfrak{u}|}} + \frac{1}{2} \left(\frac{1}{N} \prod_{j=1}^s{\left( \beta_j + \gamma_j S_N \right)}\right.\nonumber\\
									& \quad + \frac{1}{N}\sum_{p=0}^{m-1}{b^{m-p-1}(b-1) \prod_{\substack{j=1 \\ w_j \geq m-p}}^{s}{\left( \beta_j + \gamma_jS_N \right)}\prod_{\substack{j=1 \\ w_j < m-p}}^{s}{\beta_j}} - \left.\prod_{j=1}^s{\beta_j} \right).
\end{align}
\end{proof1}

\section{The reduced CBC construction}\label{sectionReducedCBC}
In this section we give a component-by-component construction for the generating vector and an upper bound for the weighted star discrepancy of the corresponding lattice rule.
\begin{alg}\label{reducedCbc}
Let $N = b^m$ and $(w_j)_{j \geq 1}$ be as above and construct $\mathbf{z} = \left( b^{w_1}z_1 , \dotsc , b^{w_s} z_s \right) \in \mathcal{Z}_{N,\mathbf{w}}^s$ as follows:
\begin{enumerate}
	\item Set $z_1 = 1$.
	\item For $d \in [s-1]$ assume $z_1, \dotsc , z_d$ to be already found. Choose $z_{d+1} \in \mathcal{Z}_{N,w_{d+1}}$ such that
	$$R_{N,\bsgamma}^{d+1}{(b^{w_1}z_1, \dotsc , b^{w_d}z_d, b^{w_{d+1}}z)}$$
	is minimized as a function of $z$.
	\item Increase $d$ by 1 and repeat the second step until $\mathbf{z} = \left( b^{w_1}z_1 , \dotsc , b^{w_s} z_s \right)$ is found.
\end{enumerate}
\end{alg}
In the algorithm above the search space is reduced for each coordinate of $\mathbf{z}$ according to its importance. This results in a considerable reduction of the construction cost as we will see in Section \ref{sectionReducedFastCBC}. This is why we call this algorithm a reduced CBC-algorithm.

The following theorem gives an upper bound for the figure of merit $R_{N, \bsgamma}^d$ of lattice point sets with generating vectors obtained from the algorithm above.
\begin{thm}\label{thm:cbc}
Let $\mathbf{z} = \left( b^{w_1}z_1 , \dotsc , b^{w_s} z_s \right)$ be constructed according to Algorithm \ref{reducedCbc}. Then for every $d \in [s]$,
\begin{equation}\label{eq:boundCbc}
R_{N,\bsgamma}^{d}{\left( b^{w_1}z_1 , \dotsc , b^{w_d} z_d \right)} \leq \frac{1}{N} \prod_{j=1}^d{\left( \beta_j + \left( 1 + 2 b^{\min{\{ w_j, m \}}}\right) \gamma_j S_N  \right)}.
\end{equation}
\end{thm}
\begin{cor}
Let $N = b^m$ and $(w_j)_{j \geq 1}$ be as above and let $\mathbf{z} = \left( b^{w_1}z_1 , \dotsc , b^{w_s} z_s \right) \in \mathcal{Z}_{N,\mathbf{w}}^s$ be constructed with Algorithm \ref{reducedCbc}. Then the corresponding lattice rule has a weighted star discrepancy
\begin{align*}
D_{N,\bsgamma}^*(\mathbf{z}) \leq \sum_{\mathfrak{u} \subseteq [s]}{\gamma_{\mathfrak{u}}\left( 1- \left( 1 - \frac{1}{N} \right) \right)^{|\mathfrak{u}|}} + \frac{1}{2N} \prod_{j=1}^{s}{\left( \beta_j + \left( 1 + 2 b^{\min{\{ w_j, m \}}}\right) \gamma_j S_N  \right)}.
\end{align*}
\end{cor}
\begin{proof1}
Combining \eqref{eq:estStar1}, \eqref{eq:S_N} and \prettyref{thm:cbc} we immediately obtain the result. 
\end{proof1}
To prove \prettyref{thm:cbc} we use the the following
\begin{lem}\label{lem:sumTnwjk}
Let $N = b^m$, $(w_j)_{j \geq 1}$ and $\mathcal{Z}_{N, w_j}$ be defined as above and recall from \eqref{eq:tnwjk} the notation
$$T_{N,w_j}(k) = \sum\limits_{z_j \in \mathcal{Z}_{N, w_j}}{ \snz{j} }.$$
Then
\begin{equation}
\sum_{k=1}^{N-1}{\frac{|T_{N,w_j}(k)|}{|\mathcal{Z}_{N, w_j}|}} \leq 2 b^{\min{\{ w_j, m \}} }S_N.
\end{equation}
\end{lem} 

\begin{proof1}
As before, we distinguish the two cases $w_j \geq m$ and $w_j < m$.\newline
{\bfseries Case 1:} $w_j \geq m$. Then \eqref{eq:tnwjk1} yields
\begin{equation*}
\sum_{k=1}^{N-1}{\frac{|T_{N,w_j}(k)|}{|\mathcal{Z}_{N, w_j}|}} = \sum_{k=1}^{N-1}{S_N} = (N-1)S_N \leq 2 N S_N = 2 b^{\min{\{ w_j, m \}} }S_N.
\end{equation*}
{\bfseries Case 2:} $w_j < m$. We use \eqref{eq:tnwjk1} and \eqref{eq:<m} to find
\begin{equation*}
\sum_{k=1}^{N-1}{\frac{|T_{N,w_j}(k)|}{|\mathcal{Z}_{N, w_j}|}} = \sum_{\substack{k=1 \\ b^{m-w_j} | k}}^{N-1}{\frac{|T_{N,w_j}(k)|}{|\mathcal{Z}_{N, w_j}|}} + \sum_{\substack{k=1 \\ b^{m-w_j} \nmid \, k}}^{N-1}{\frac{|T_{N,w_j}(k)|}{|\mathcal{Z}_{N, w_j}|}} = (b^{w_j}-1)S_N + b^{w_j}\sum_{r=1}^{b^{m-w_j}-1}{\frac{|T_{N,w_j}(r)|}{|\mathcal{Z}_{N, w_j}|}}.
\end{equation*}
For any $r \in \{1, \dotsc , b^{m-w_j}-1\}$ the condition $\gcd{(r, b^{m-w_j})} = b^{\nu}$ is equivalent to $b^{\nu} \mid r$ and $b^{\nu + 1} \nmid r$ simultaneously. Using this we investigate the last sum in the above equation:
\begin{align*}
\sum_{r=1}^{b^{m-w_j}-1}{\frac{|T_{N,w_j}(r)|}{|\mathcal{Z}_{N, w_j}|}} &= \frac{1}{|\mathcal{Z}_{N, w_j}|}\sum_{\nu = 0}^{m-w_j-1}{\sum_{\substack{r=1 \\ b^{\nu} \mid r \\ b^{\nu +1} \nmid \, r}}^{b^{m-w_j}-1}{|T_{N,w_j}(r)|}}.
\end{align*}
Once more with the aid of \eqref{eq:tnwjk1} this yields
\begin{align*}
\sum_{r=1}^{b^{m-w_j}-1}{\frac{|T_{N,w_j}(r)|}{|\mathcal{Z}_{N, w_j}|}} &= \frac{1}{|\mathcal{Z}_{N, w_j}|}\sum_{\nu = 0}^{m-w_j-1}{\sum_{\substack{r=1 \\ b^{\nu} \mid r \\ b^{\nu +1} \nmid \, r}}^{b^{m-w_j}-1}{\left|b^{\nu}(S_{b^{w_j + \nu}}-S_{b^{w_j + \nu + 1}}) \right|}} \\
																																				&= \frac{1}{|\mathcal{Z}_{N, w_j}|}\sum_{\nu = 0}^{m-w_j-1}{\sum_{\substack{r=1 \\ b^{\nu} \mid r \\ b^{\nu +1} \nmid \, r}}^{b^{m-w_j}-1}{b^{\nu}(S_{b^{w_j + \nu + 1}} - S_{b^{w_j + \nu}}) }}.
\end{align*}
Analogously to \eqref{eq:cardinality||} we find
\begin{align*}
\left| \left\{ r \in \left\{ 1, \dotsc , b^{m-w_j} -1 \right\} : b^{\nu} \mid r \text{ and } b^{\nu+1} \nmid r \right\} \right| = b^{m-w_j-\nu-1}(b-1)
\end{align*}
and hence
\begin{align*}
\sum_{r=1}^{b^{m-w_j}-1}{\frac{|T_{N,w_j}(r)|}{|\mathcal{Z}_{N, w_j}|}}	&= \sum_{\nu = 0}^{m-w_j-1}{(S_{b^{w_j + \nu + 1}} - S_{b^{w_j + \nu}}) } = S_N - S_{b^{w_j}}.
\end{align*}

Altogether we have
\begin{align*}
\sum_{k=1}^{N-1}{\frac{|T_{N,w_j}(k)|}{|\mathcal{Z}_{N, w_j}|}} = (b^{w_j}-1)S_N + b^{w_j}(S_N - S_{b^{w_j}}) \leq 2 b^{w_j}S_N = 2 b^{\min{\{ w_j, m \}} }S_N
\end{align*}
and the proof is complete.
\end{proof1}

With the aid of \prettyref{lem:sumTnwjk} we are able to prove \prettyref{thm:cbc} using induction on $d$.
\begin{proof1}
According to Algorithm \ref{reducedCbc} we set $z_1=1$ in Step 1. We have to show that
$$R_{N,\bsgamma}^1(b^{w_1}) \leq \frac{1}{N}\left( \beta_1 + \left( 1 + 2 b^{\min\{ w_1, m \}}\right) \gamma_1 S_N \right).$$
With \eqref{eq:endSquared1} we have
\begin{equation*}
R_{N,\bsgamma}^1(b^{w_1}) = \frac{1}{N} \sum_{k=0}^{N-1}{\left( \beta_1 + \gamma_1 \snb{w_1} \right)} - \beta_1 = \frac{1}{N} \sum_{k=0}^{N-1}{\gamma_1 \snb{w_1}}.
\end{equation*}
Again, we consider the two cases $w_1 \geq m$ and $w_1 < m$ separately.\newline
{\bfseries Case 1:} $w_1 \geq m$. Then
\begin{equation*}
\begin{split}
R_{N,\bsgamma}^1(b^{w_1})	&= \frac{1}{N} \sum_{k=0}^{N-1}{\gamma_1 \snbn{w_1-m}} = \frac{1}{N}\gamma_1 N S_N \leq \frac{1}{N}(1 + \gamma_1 + 2 N \gamma_1 S_N) \\
										&= \frac{1}{N}\left(\beta_1 + 2 b^{\min\{ w_1, m \}} \gamma_1 S_N\right) \leq \frac{1}{N}\left( \beta_1 + \left( 1 + 2 b^{\min\{ w_1, m \}}\right) \gamma_1 S_N \right)
\end{split}
\end{equation*}
which is the desired result.\newline
{\bfseries Case 2:} $w_1 < m$. After interchanging the two sums, once more, we split up the inner sum as follows:
\begin{equation*}
\begin{split}
R_{N,\bsgamma}^1(b^{w_1})	&= \frac{\gamma_1}{N} \sn \sum_{k=0}^{N-1}{\mathrm{e}^{2 \pi \ii h k/b^{m-w_1}}} \\
										&= \frac{\gamma_1}{N} \sum_{\substack{-\frac{N}{2} < h \leq \frac{N}{2} \\ h \neq 0 \\ b^{m-w_1} \mid \, h}}{\frac{1}{|h|} \sum_{k=0}^{N-1}{\mathrm{e}^{2 \pi \ii h k/b^{m-w_1}}}} + \frac{\gamma_1}{N} \sum_{\substack{-\frac{N}{2} < h \leq \frac{N}{2} \\ h \neq 0 \\ b^{m-w_1} \nmid \, h}}{\frac{1}{|h|} \sum_{k=0}^{N-1}{\mathrm{e}^{2 \pi \ii h k/b^{m-w_1}}}}.
\end{split}
\end{equation*}
Now we are able to compute the inner sums. The first one sums to $N$, whereas the second one equals zero which can immediately be seen by applying the formula for finite geometric series. Thus
\begin{equation*}
\begin{split}
R_{N,\bsgamma}^1(b^{w_1})	&= \gamma_1 \sum_{\substack{-\frac{N}{2} < h \leq \frac{N}{2} \\ h \neq 0 \\ b^{m-w_1} \mid \, h}}{\frac{1}{|h|}}.
\end{split}
\end{equation*}
We use \eqref{eq:NiedSn} to find
\begin{equation*}
\begin{split}
R_{N,\bsgamma}^1(b^{w_1})	&= \gamma_1 \frac{1}{b^{m-w_1}}S_{\frac{N}{b^{m-w_1}}} = \frac{\gamma_1}{N}b^{w_1}S_{b^{w_1}} \leq \frac{\gamma_1}{N}b^{w_1}S_{N} \leq \frac{1}{N}\left( \beta_1 + 2 b^{w_1} \gamma_1 S_N \right) \\
										&\leq \frac{1}{N}\left( \beta_1 + \left( 1 + 2 b^{\min\{ w_1, m \}}\right) \gamma_1 S_N \right),
\end{split}
\end{equation*}
as claimed.

Let $d \in [s-1]$ and assume that we have some $\mathbf{z} \in \mathcal{Z}_{N,\mathbf{w}}^d$, such that
$$R_{N,\bsgamma}^d{\left( b^{w_1}z_1 , \dotsc , b^{w_d} z_d \right)} \leq \frac{1}{N} \prod_{j=1}^d{\left( \beta_j + \left( 1 + 2 b^{\min{\{ w_j, m \}}}\right) \gamma_j S_N  \right)}.$$

We have to prove the existence of a $z_{d+1} \in \mathcal{Z}_{N,w_{d+1}}$ with
$$R_{N,\bsgamma}^{d+1}{\left( b^{w_1}z_1 , \dotsc , b^{w_d} z_d, b^{w_{d+1}} z_{d+1} \right)} \leq \frac{1}{N} \prod_{j=1}^{d+1}{\left( \beta_j + \left( 1 + 2 b^{\min{\{ w_j, m \}}}\right) \gamma_j S_N  \right)}.$$
Using again \eqref{eq:endSquared1} we have for any $z_{d+1} \in \mathcal{Z}_{N,w_{d+1}}$ that
\begin{equation}\label{eq:endSquared2}
\begin{split}
R	&_{N,\bsgamma}^{d+1}{\left( b^{w_1}z_1 , \dotsc , b^{w_d} z_d, b^{w_{d+1}} z_{d+1} \right)} \\
	&= \frac{1}{N}\sum_{k=0}^{N-1}{\prod_{j=1}^d{\left( \beta_j + \gamma_j \snz{j} \right)}\left( \beta_{d+1} + \gamma_{d+1} \snz{d+1} \right)} - \beta_{d+1} \prod_{j=1}^d{\beta_j}\\
	&= \beta_{d+1} R_{N,\bsgamma}^d{\left( b^{w_1}z_1 , \dotsc , b^{w_d} z_d \right)} \\
	&\quad + \frac{\gamma_{d+1}}{N}\sum_{k=0}^{N-1}{\prod_{j=1}^d{\left( \beta_j + \gamma_j \snz{j} \right)}} \snz{d+1}\\
	&= \beta_{d+1} R_{N,\bsgamma}^d{\left( b^{w_1}z_1 , \dotsc , b^{w_d} z_d \right)} + \frac{\gamma_{d+1} S_N}{N} \prod_{j=1}^d{\left( \beta_j + \gamma_j S_N \right)} \\
	&\quad + \frac{\gamma_{d+1}}{N}\sum_{k=1}^{N-1}{ \snz{d+1} \prod_{j=1}^d{\left( \beta_j + \gamma_j \snz{j} \right)}}.				
\end{split}
\end{equation}

Next we consider the arithmetic mean of $R_{N,\bsgamma}^d{\left( b^{w_1}z_1 , \dotsc , b^{w_d} z_d, b^{w_{d+1}} z \right)}$ over all $z \in \mathcal{Z}_{N,w_{d+1}}$. As only the third summand in \eqref{eq:endSquared2} depends on the $(d+1)$-st coordinate it suffices to investigate the mean of this summand. Clearly, if we have some upper bound for the mean over all $z \in \mathcal{Z}_{N,w_{d+1}}$, there exists a $z_{d+1} \in \mathcal{Z}_{N,w_{d+1}}$ which satisfies this bound. Thus we study
\begin{align*}
\frac{1}{|\mathcal{Z}_{N,w_{d+1}}|}	&\sum_{z \in \mathcal{Z}_{N,w_{d+1}}}{\frac{\gamma_{d+1}}{N}\sum_{k=1}^{N-1}{ \snzo{d+1} \prod_{j=1}^d{\left( \beta_j + \gamma_j \snz{j} \right)}}}.
\end{align*}
We bound this term by its absolute value
\begin{align*}
&\left| \frac{1}{|\mathcal{Z}_{N,w_{d+1}}|} \sum_{z \in \mathcal{Z}_{N,w_{d+1}}}\frac{\gamma_{d+1}}{N}\sum_{k=1}^{N-1}\snzo{d+1} \prod_{j=1}^d{\left( \beta_j + \gamma_j \snz{j} \right)} \right| \\
																		&\leq \frac{\gamma_{d+1}}{N} \sum_{k=1}^{N-1}{\frac{1}{|\mathcal{Z}_{N,w_{d+1}}|} \left| \sum_{z \in \mathcal{Z}_{N,w_{d+1}}}{ \snz{d+1} } \right|}\\
																		&\quad \times \prod_{j=1}^d{\left( \beta_j + \gamma_j \snzabs{j} \right)} \\
																		&\leq \frac{\gamma_{d+1}}{N}\sum_{k=1}^{N-1}{\frac{|T_{N, w_{d+1}}(k)|}{|\mathcal{Z}_{N,w_{d+1}}|}}\prod_{j=1}^d{\left( \beta_j + \gamma_j S_N \right)} \\
																		&\leq \frac{\gamma_{d+1}}{N}2 b^{\min\{ w_{d+1}, m \}}S_N \prod_{j=1}^d{\left( \beta_j + \gamma_j S_N \right)},
\end{align*}
where the last estimate stems from application of \prettyref{lem:sumTnwjk}. Combining this with \eqref{eq:endSquared2} we have shown the existence of a $z_{d+1} \in \mathcal{Z}_{N,w_{d+1}}$ such that
\begin{equation*}
\begin{split}
R	&_{N,\bsgamma}^{d+1}{\left( b^{w_1}z_1 , \dotsc , b^{w_d} z_d, b^{w_{d+1}} z_{d+1} \right)} \\
	&\leq \beta_{d+1} R_{N,\bsgamma}^d{\left( b^{w_1}z_1 , \dotsc , b^{w_d} z_d \right)} + \frac{\gamma_{d+1} S_N}{N} \prod_{j=1}^d{\left( \beta_j + \gamma_j S_N \right)} + \frac{\gamma_{d+1}}{N}2 b^{\min\{ w_{d+1}, m \}} S_N \prod_{j=1}^d{\left( \beta_j + \gamma_j S_N \right)}.
\end{split}
\end{equation*}
We use the induction hypothesis to find
\begin{equation*}
\begin{split}
R	&_{N,\bsgamma}^{d+1}{\left( b^{w_1}z_1 , \dotsc , b^{w_d} z_d, b^{w_{d+1}} z_{d+1} \right)} \\
	&\leq \frac{\beta_{d+1}}{N} \prod_{j=1}^d{\left( \beta_j + \left( 1 + 2 b^{\min{\{ w_j, m \}}}\right) \gamma_j S_N  \right)} + \frac{\gamma_{d+1} S_N}{N} \prod_{j=1}^d{\left( \beta_j + \gamma_j S_N \right)}\left(1 + 2 b^{\min\{ w_{d+1}, m \}} \right) \\
	&\leq \left( \beta_{d+1} + \left(1 + 2 b^{\min\{ w_{d+1}, m \}} \right) \gamma_{d+1} S_N \right) \frac{1}{N} \prod_{j=1}^d{\left( \beta_j + \left( 1 + 2 b^{\min{\{ w_j, m \}}}\right) \gamma_j S_N  \right)} \\
	&= \frac{1}{N} \prod_{j=1}^{d+1}{\left( \beta_j + \left( 1 + 2 b^{\min{\{ w_j, m \}}}\right) \gamma_j S_N  \right)}
\end{split}
\end{equation*}
which completes the proof.
\end{proof1}

\section{The reduced fast CBC construction}\label{sectionReducedFastCBC}
By now we have seen how we can construct a generating vector of a lattice point set with low weighted star discrepancy with a reduced CBC construction as in the previous section. Now we study the construction cost of this algorithm. In fact the CBC algorithm can be made faster to construct generating vectors for relatively large $N$ and $s$. To show this we follow closely \cite{DKLP} and \cite{LP}.

Let $d \in [s-1]$ and assume that we have already found $(b^{w_1}z_1, \dotsc , b^{w_d}z_d)$. Then we have (cf. \eqref{eq:endSquared1})
$$R_{N,\bsgamma}^d(b^{w_1}z_1, \dotsc , b^{w_d}z_d) = \frac{1}{N}\sum_{k=0}^{N-1}{\prod_{j=1}^d{\left( \beta_j + \gamma_j \snz{j} \right)}} - \prod_{j=1}^d{\beta_j}.$$
Define $r(h) = \max{\{ 1, |h| \}}$. Then
\begin{align*}
\beta_j + \gamma_j \snz{j}	&= \beta_j + \gamma_j\left( \sum_{-\frac{N}{2} < h \leq \frac{N}{2}}{\frac{\mathrm{e}^{2 \pi \ii h k b^{w_j}z_j/N}}{r(h)}} - 1 \right)\\
														&= 1 + \gamma_j \sum_{-\frac{N}{2} < h \leq \frac{N}{2}}{\frac{\mathrm{e}^{2 \pi \ii h k b^{w_j}z_j/N}}{r(h)}}.
\end{align*}
Hence we have
\begin{equation}\label{eq:fastCBC}
\begin{split}
R_{N,\bsgamma}^d(b^{w_1}z_1, \dotsc , b^{w_d}z_d) &= \frac{1}{N}\sum_{k=0}^{N-1}\prod_{j=1}^d{\left(1 + \gamma_j \sum_{-\frac{N}{2} < h \leq \frac{N}{2}}{ \frac{\mathrm{e}^{2 \pi \ii h k b^{w_j}z_j/N}}{r(h)}}\right)} - \prod_{j=1}^d{\beta_j}\\
																						&= \frac{1}{N}\sum_{k=0}^{N-1}\eta_d(k) - \prod_{j=1}^d{\beta_j},
\end{split}
\end{equation}
where we have defined
$$\eta_d(k) = \prod_{j=1}^d{\left(1 + \gamma_j \phi\left(\frac{k b^{w_j}z_j}{N}\right)\right)}$$
and
$$\phi(x) = \sum_{-\frac{N}{2} < h \leq \frac{N}{2}}{ \frac{\mathrm{e}^{2 \pi \ii h x}}{r(h)}}.$$

However, this is exactly the situation as dealt with in \cite[Section~4.2]{LP}. Thus we know that $\phi\left(\frac{k b^{w_j}z_j}{N}\right)$ takes on at most $N$ different values, namely
$$\phi(0), \phi\left( \frac{1}{N} \right), \dotsc , \phi\left( \frac{N-1}{N} \right),$$ 
which can be computed in $O(N \log N)$ operations and stored in a memory space of size $O(N)$, as demonstrated in \cite{LP}.

Next we investigate one actual step of the CBC construction. Assuming that we have already found $(b^{w_1}z_1, \dotsc , b^{w_d}z_d) \in \mathcal{Z}_{N,\mathbf{w}}^d$ we have to minimize $R_{N,\bsgamma}^{d+1}(b^{w_1}z_1, \dotsc , b^{w_d}z_d,b^{w_{d+1}}z)$ as a function of $z \in \mathcal{Z}_{N,w_{d+1}}$ to find $z_{d+1} \in \mathcal{Z}_{N,w_{d+1}}$. For $w_{d+1} \geq m$ we just set $z_{d+1} = 1$ and are done. Therefore let $w_{d+1} < m$. Considering \eqref{eq:fastCBC} we have
\begin{align*}
R_{N,\bsgamma}^{d+1}(b^{w_1}z_1, \dotsc , b^{w_d}z_d,b^{w_{d+1}}z_{d+1})	&= \frac{1}{N}\sum_{k=0}^{N-1}\eta_{d+1}(k) - \prod_{j=1}^{d+1}{\beta_j}\\
																																					&= \frac{1}{N}\sum_{k=0}^{N-1}{\eta_{d}(k)\left(1 + \gamma_{d+1} \phi\left(\frac{k b^{w_{d+1}}z_{d+1}}{N}\right)\right)} - \prod_{j=1}^{d+1}{\beta_j}\\
																																					&= \frac{1}{N}\sum_{k=0}^{N-1}{\eta_{d}(k)\left(1 + \gamma_{d+1} \phi\left(\left\{ \frac{k b^{w_{d+1}}z_{d+1}}{N}\right\}\right)\right)} - \prod_{j=1}^{d+1}{\beta_j}.
\end{align*}
It is obviously enough to minimize $\frac{1}{N}\sum_{k=0}^{N-1}{\eta_{d}(k)\phi\left(\left\{ \frac{k b^{w_{d+1}}z_{d+1}}{N}\right\}\right)}$. To do this we proceed analogously to \cite{DKLP}. We define the matrix
$$A = \left( \phi{\left( \left\{ \frac{k b^{w_{d+1}}z_{d+1}}{N}\right\}\right)} \right)_{z_{d+1} \in \mathcal{Z}_{N,w_j}, k \in \{ 0, \dotsc , N-1 \}}$$
and
$$\bseta_{d} = \left( \eta_{d}(0), \eta_d(1), \dotsc , \eta_d(N-1) \right)^{\top}$$
and find that
$$\frac{1}{N}\sum_{k=0}^{N-1}{\eta_{d}(k)\left(\phi\left(\left\{ \frac{k b^{w_{d+1}}z_{d+1}}{N}\right\}\right)\right)} = A \bseta_{d}.$$
We can display the matrix $A$ as 
$$A = (\Omega^{(m-w_{d+1})}, \dotsc , \Omega^{(m-w_{d+1})}),$$
with
$$\Omega^{(l)} = \left(\phi\left(\left\{ \frac{k z_{d+1}}{b^l}\right\}\right)\right)_{z_{d+1} \in \mathcal{Z}_{b^l,0}, k \in \{ 0, \dotsc , b^l - 1 \}}.$$
Again analogously to \cite{DKLP} we obtain the following reduced fast CBC algorithm.
\begin{alg}\label{reducedFastCBC}
\begin{enumerate}[a)]
\item Compute $\phi\left( \frac{r}{N} \right)$ for all $r = 0, \dotsc , N-1.$
\item Set $\eta_1(k) = 1 + \gamma_1\phi\left( \left\{ \frac{k b^{w_1} z_1}{N} \right\} \right)$ for $k = 0, \dotsc , N-1.$
\item Set $z_1 = 1$. Set $d = 2$ and recall that we have defined $t = \max\{j : w_j < m\}$. While $d \leq \min\{ s, t\}$,
	\begin{enumerate}[1.]
	\item partition $\bseta_{d-1}$ into $b^{w_d}$ vectors $\bseta_{d-1}^{(1)}, \dotsc , \bseta_{d-1}^{(b^{w_d})}$ of length $b^{m-w_d}$ and let $\bseta' = \bseta_{d-1}^{(1)} + \dotsc + \bseta_{d-1}^{(b^{w_d})}$ denote their sum,
	\item let $T_d(z) = \Omega^{(m - w_d)}\bseta'$,
	\item let $z_d = \mathrm{arg}\,\mathrm{min}_z{T_d(z)}$,
	\item let $\eta_d(k) = \eta_{d-1}\left(1 + \gamma_d\phi\left( \left\{ \frac{k b^{w_d} z_d}{N} \right\} \right)\right)$ for $k = 0, \dotsc , N-1$,
	\item increase $d$ by 1.
	\end{enumerate}
If $s > t$, then set $z_t = z_{t+1} = \dotsc = z_s = 1$. Then we have
$$R_{N,\bsgamma}^s\left( b^{w_1}z_1, \dotsc , b^{w_s}z_s \right) = \frac{1}{N}\sum_{k=0}^{N-1}{\eta_s(k)}.$$
\end{enumerate}
\end{alg}
Using \cite{DKLP, LP, NC1, NC2} we find that Algorithm \ref{reducedFastCBC} has a construction cost of
$$O\left( N \log N + \min\{s,t\} N + N \sum_{d=1}^{\min\{s,t\}}{(m-w_d)b^{-w_d}} \right)$$
operations, in comparison to $O(sN\log N)$ operations for the standard CBC algorithm used for example in \cite{SJ1}.

\section{Conditions for strong polynomial tractability}\label{sectionTractability}
Let $\mathbf{z} = \left( b^{w_1}z_1 , \dotsc , b^{w_s} z_s \right) \in \mathcal{Z}_{N,\mathbf{w}}^s$ be constructed with Algorithm \ref{reducedCbc} or \ref{reducedFastCBC} and consider the corresponding lattice rule. We are interested in conditions for tractability of the weighted star discrepancy of such lattice point sets. 
From \eqref{eq:estStar1} and \eqref{eq:endSquared} we know
$$D_{N, \bsgamma}^*(\mathbf{z}) \leq \sum_{\mathfrak{u} \subseteq [s]}{\gamma_{\mathfrak{u}}\left( 1- \left( 1 - \frac{1}{N} \right) \right)^{|\mathfrak{u}|}} + \frac{1}{2} R_{N,\bsgamma}^s(\mathbf{z}).$$

For now, let us assume that the $\gamma_jb^{w_j}$'s are summable, i.e. $\sum\limits_{j=1}^{\infty}{\gamma_jb^{w_j}} < \infty.$ Similar to Joe and Sinescu in \cite{J} and \cite{SJ1}, we see that in this case \eqref{eq:estStar1} implies
$$D_{N, \bsgamma}^*(\mathbf{z}) \leq O\left( \frac{1}{N} \right) + \frac{1}{2} R_{N,\bsgamma}^s(\mathbf{z}),$$
where the implied constant in the $O$-notation is independent of $s$.

Recall that we have defined the information complexity as
$$N^*(\varepsilon,s) := \min\left\{ N \in \mathbbm{N} : D_{N, \bsgamma}^*(\mathbf{z}) \leq \varepsilon \right\}.$$
If $\sum\limits_{j=1}^{\infty}{\gamma_jb^{w_j}} < \infty$, it is easy to show that it is equivalent to consider the standard notions of tractability with respect to $N^*(\varepsilon,s)$ or with respect to $\min\left\{ N \in \mathbbm{N} : R_{N,\bsgamma}^s(\mathbf{z}) \leq \varepsilon \right\}$.

\prettyref{thm:cbc} yields
\begin{equation*}
R_{N,\bsgamma}^s(\mathbf{z}) \leq \frac{1}{N} \prod_{j=1}^s{\left( \beta_j + \left( 1 + 2 b^{\min{\{ w_j, m \}}}\right) \gamma_j S_N  \right)}.
\end{equation*}
We study the right-hand side of the latter inequality.
\begin{equation}\label{eq:RHS}
\begin{split}
\frac{1}{N} \prod_{j=1}^s{\left( \beta_j + \left( 1 + 2 b^{\min{\{ w_j, m \}}}\right) \gamma_j S_N  \right)}	&\leq \frac{1}{N} \prod_{j=1}^s{\left( \beta_j + \left( 1 + 2 b^{\min{\{ w_j, m \}}}\right) \gamma_j 2\left( \log{\left\lfloor \frac{N}{2} \right\rfloor} + 1 \right)  \right)}\\
																																																							&\leq \frac{1}{N} \prod_{j=1}^s{\left( \beta_j + \left( 1 + 2 b^{\min{\{ w_j, m \}}}\right) \gamma_j 4\log N  \right)}\\
																																																							&= \frac{1}{N} \prod_{j=1}^s{\left( 1 + \gamma_j\left(1 + 4 \left(1 + 2 b^{\min{\{ w_j, m \}}}\right)\log N\right)  \right)},
\end{split}
\end{equation}
where we have used 
$$S_N = \sn \leq 2 \sum_{h=1}^{\left\lfloor \frac{N}{2} \right\rfloor}{\frac{1}{h}} \leq 2 \log{\left\lfloor \frac{N}{2} \right\rfloor} + 2 \leq 4 \log N.$$
The second to last inequality is a well-known estimate for partial sums of the harmonic series.

Now we have
\begin{align*}
\frac{1}{N} \prod_{j=1}^s{\left( \beta_j + \left( 1 + 2 b^{\min{\{ w_j, m \}}}\right) \gamma_j S_N  \right)}	&\leq \frac{1}{N} \prod_{j=1}^s{\left( 1 + \gamma_j\left(1 + 4 \left(1 + 2 b^{w_j}\right)\log N\right)  \right)}\\
																																																							&\leq \frac{1}{N} \prod_{j=1}^s{\left( 1 + 13 \gamma_j b^{w_j} \log N  \right)}.
\end{align*}
Define $\sigma_d := 13 \sum\limits_{j=d+1}^{\infty}{\gamma_jb^{w_j}}$ for $d\geq 0$. From \cite[p. 222]{DP} or \cite[Lemma~3]{HN} we know that
$$\prod_{j=1}^s{\left( 1 + 13 \gamma_jb^{w_j} \log N  \right)} \leq \left( 1 + \sigma_d^{-1} \right)^d N^{(\sigma_0 + 1)\sigma_d}.$$
For $0 < \delta < 1$ choose $d$ large enough such that $\sigma_d \leq \frac{\delta}{\sigma_0 + 1}$. Then
$$\prod_{j=1}^s{\left( 1 + 13 \gamma_jb^{w_j} \log N  \right)} \leq c_{\bsgamma, \delta}N^{\delta},$$
where $c_{\bsgamma, \delta}$ is independent of $s$ and $N$.
Thus we have
$$R_{N,\bsgamma}^s(\mathbf{z}) \leq c_{\bsgamma, \delta}N^{\delta - 1}.$$
We obtain $c_{\bsgamma, \delta}N^{\delta - 1} \leq \varepsilon$ and thus $R_{N,\bsgamma}^s(\mathbf{z}) \leq \varepsilon$ if $N \geq (c_{\bsgamma, \delta} \varepsilon^{-1})^{\frac{1}{1-\delta}}$. Hence, if the $\gamma_j b^{w_j}$'s are summable we always achieve strong polynomial tractability.
\begin{rem}
Whether the conditions on the $\gamma_j$'s and $w_j$'s can be mitigated while at least polynomial tractability still holds remains an unresolved problem.
\end{rem}

\section*{Acknowledgements}{\small
The authors would like to thank Peter Kritzer and Friedrich Pillichshammer for their comments and suggestions.}

\noindent
{\bfseries Addresses:}

R. Kritzinger, H. Laimer, Department of Financial Mathematics and Applied Number Theory, Johannes Kepler University Linz, Altenbergerstr. 69, 4040 Linz, Austria.\newline
e-mail: ralph.kritzinger@jku.at, helene.laimer@jku.at
\end{document}